# Non-trivial Local Attractors of a Three-dimensional Dynamical System

Keying Guan

Science College, Beijing Jiaotong University,

Beijing, China, 100044

Email: keying.guan@gmail.com

Abstract: Based on both qualitative method and numerical tests for a series of particular cases in the parameter region, $a = 1$, $0 < b < 1$, it is shown that the three-dimensional system (2) may have a series of interesting phenomena on the non-trivial local attractors, such as the "faint attractor" (this term is suggested by the author), the local attractor with complex structure, twin spatial limit closed orbits, the bifurcation of rotation numbers, and the spatial limit cycle, etc.. The system (2) is a very rich source in the study of dynamical system theory.

In his research about the diffusion problem of small amplitude disturbance on a solid surface, professor R. Mark Bradley, posted a question: "How can I solve a^2 y = y^3 − y'''− y' ?" on the website of ResearchGate:

https://www.researchgate.net/post/How_can_I_solve_a2_y_y3-y-y?cp=re65_x_p2&ch=reg&loginT=v4u08vKuAj1OFM-Xv1OU0FSa98kmnppZC555yxiDD-0%2C&pli=1#view=527b088fd2fd6427258b4868

The equation in his problem is a third order ordinary differential equation.

$$y''' + y' + a^2 y - y^3 = 0 \qquad (1)$$

In the discussion about his question, I noticed that (1) is corresponding to a 3-dimensional dynamical system with three interesting equilibrium, but the system does not have any attractor for it is conservative. I found that by adding a new term, $b\, y''$, to (1), it then will correspond a new three-dimensional dynamical system

$$\begin{cases} \dfrac{dx}{dt} = P(x,y,z) = y \\ \dfrac{dy}{dt} = Q(x,y,z) = z \\ \dfrac{dz}{dx} = R(x,y,z) = x^3 - a^2\, x - y - b\, z \end{cases} \qquad (2)$$

where the parameters $a$ and $b$ are both positive numbers. This paper will uncover a series of interesting facts about its attractor. The main knowledge on the dynamical systems used in this paper can be found in the famous book of John Guckenheimer and Philip Holmes, Nonlinear Oscillations, Dynamical Systems, and Bifurcations of Vector Fields.

## § 1. Basic Properties of the System (2)

The system (2) is an autonomous polynomial system with three equilibrium points: $p_0 = (0,0,0)$, $p_1 = (-a, 0,0)$ and $p_2 = (a, 0,0)$. When $b = 0$, the system is equivalent to (1). However, when $b > 0$, the system (2) is dissipative, since

$$\frac{\partial P}{\partial x} + \frac{\partial Q}{\partial y} + \frac{\partial R}{\partial z} = -b < 0 \tag{3}$$

The system (2) is similar to the famous Lorenz equation

$$\begin{cases} \frac{dx}{dt} = \sigma(y - x) \\ \frac{dy}{dt} = \rho x - y - xz \\ \frac{dz}{dt} = -\beta z + xy \end{cases} \tag{4}$$

in the follow sense:

a. They have similar symmetries, i.e., (2) is symmetric about the origin, for it is invariant when $x \to -x$, $y \to -y$ and $z \to -z$, while (4) is symmetric about the z-axis, for it is invariant when $x \to -x$, $y \to -y$ and $z \to z$.

b. Both have three equilibrium points. The three equilibrium points of (4) are $q_0 = (0,0,0)$, $q_1 = (-\sqrt{\beta(\rho - 1)}, -\sqrt{\beta(\rho - 1)}, \rho - 1)$, $q_2 = (-\sqrt{\beta(\rho - 1)}, -\sqrt{\beta(\rho - 1)}, \rho - 1)$.

c. Both are dissipative systems.

Because of the dissipative, the system (2) may have also some attractors. Clearly, since (2) is still different from (4) in several aspects, the attractors should have some new properties different from the Lorenz's attractor.

For the system (2), this paper discusses only some particular cases in the parameter region, $a = 1$, $0 \leq b < 1$.

In this region, all of the three equilibrium points are saddle-focus type. Concretely, the equilibrium point $p_0$ has a one-dimensional stable manifold $M_{1d}^s(p_0)$ and a two-dimensional unstable manifold $M_{2d}^u(p_0)$ that is corresponding to a pair of conjugate eigenvalues with positive real part, while both equilibrium points $p_1$ and $p_2$ have, respectively, a one-dimensional unstable manifold $M_{1d}^u(p_i)$ and a two-dimensional stable manifold $M_{2d}^s(p_i)$ (i=1,2) that corresponds to a pair of conjugate eigenvalues with negative real part.

In order to see these manifolds intuitively, based on numerical method, this paper uses differently colored integral curves to represent them, i.e., using light green integral curves to expand the two-dimensional stable manifolds, $M_{2d}^s(p_1)$ and $M_{2d}^s(p_2)$, using the purple integral curves to expand the two-dimensional unstable manifold $M_{2d}^u(p_0)$, using dark green to color the one-dimensional unstable manifolds, $M_{1d}^u(p_1)$ and $M_{1d}^u(p_2)$, and using red to color the one-dimensional stable manifold $M_{1d}^s(p_0)$. The three equilibrium points are colored black. An integral curve, which does not belong to any above-mentioned manifolds, is colored blue.

As a concrete example, in the case $a = 1$, $b = 0.3$, the Figures 1, 2, 3 show these manifolds of the three equilibrium points respectively, and Figure 4 puts them together by their real positions.

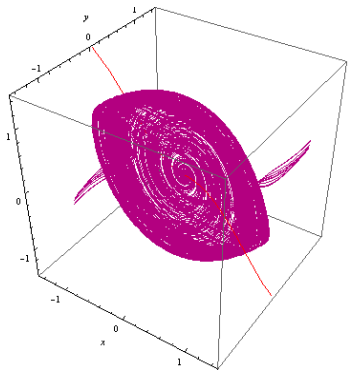
Figure 1. $M_{1d}^s(p_0)$ and $M_{2d}^u(p_0)$

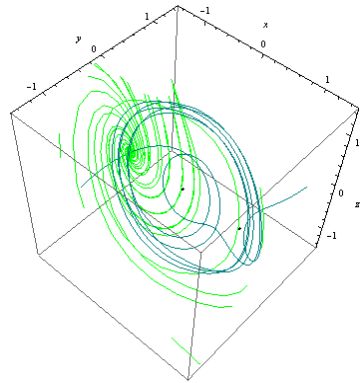
Figure 2. $M_{1d}^u(p_1)$ and $M_{2d}^s(p_1)$

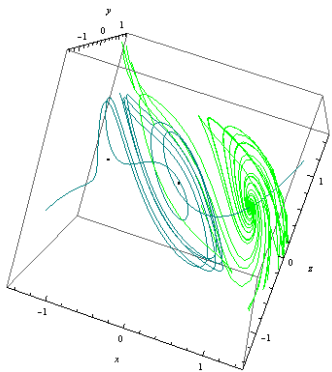
Figure 3. $M_{1d}^u(p_2)$ and $M_{2d}^s(p_2)$

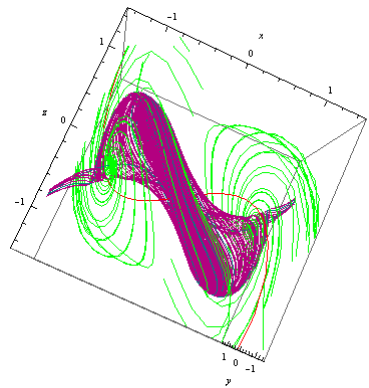
Figure 4. The six manifolds

Largely, the properties of these manifolds and the relations between them, determine the main behaviors of the system (2).

By a simple, local qualitative investigation, it is easy to see that, the x-axes pass the three two-dimensional manifolds, $M_{2d}^u(p_0)$ and $M_{2d}^s(p_i)$ (i=1,2) transversely through their corresponding equilibrium points, and the integral curves on these manifolds rotate around x-axes with the same left-hand rotational direction (relative to the positive direction of x-axis), or in other words, all of the rotations are counterclockwise if you watch them from the viewing point (5,0,0). On $M_{2d}^u(p_0)$, all integral curves rotate spirally away from the equilibrium point (in other words, all integral curves on it approach spirally to $p_0$ as $t \to -\infty$), while on $M_{2d}^s(p_i)$, all integral curves approach spirally to the corresponding equilibrium point $p_i$ (i=1,2) as $t \to +\infty$.

For a saddle-focus equilibrium point, its one-dimensional stable (or unstable) manifold must be separated into two branches by its two-dimentional unstable (or stable) manifold. It is easy to see that the two branches of $M_{1d}^s(p_0)$ are symmetric about $p_0$.

The two branches of $M_{1d}^u(p_1)$ are not symmetric about $p_1$. One branch, noted by $M_{1d}^u(p_1)_1$, goes first towards the region near the origin, and the other branch, noted by $M_{1d}^u(p_1)_2$, clearly goes far from the origin. $M_{1d}^u(p_1)_1$ is more curved than $M_{1d}^u(p_1)_2$.

Likewise, the two branches of $M_{1d}^u(p_2)$ are not symmetric about $p_2$. One branch, noted by $M_{1d}^u(p_2)_1$, goes first toward the region near the origin, and the other branch, noted by $M_{1d}^u(p_2)_2$, clearly goes far from the origin. $M_{1d}^u(p_2)_1$ is more curved than $M_{1d}^u(p_2)_2$.

However, $M_{1d}^u(p_1)$ and $M_{1d}^u(p_2)$ are symmetric about the origin, and $M_{2d}^s(p_1)$ and $M_{2d}^s(p_2)$ are also symmetric about the origin.

It is clear from the relevant definitions that there cannot be an intersection between any two different stable manifolds, and also that there cannot be an intersection between any two different unstable manifolds.

In some particular cases, a stable manifold and an unstable manifold may intersect through the homoclinic orbit, if these two manifolds belong to the same equilibrium point, or through the heteroclinic orbit, if these two manifolds belong to different equilibrium points. The system is not structurally stable under these particular cases. Though it is very important to find out these particular cases in the study of the dynamical system theory, it is usually very difficult to find out analytically the exact parameters under which these structurally unstable phenomena just happen. Therefore, this paper shall avoid the very fine work to find out the exact condition of these particular cases. However, we will guess at the possible regions in which some of these particular cases may happen.

In this paper, when the closure of the unstable manifold $M_{2d}^u(p_0)$ is a bounded set, the local attractor(s) of the system (2) will be given in the following way, that is, to establish first a class $A_{class}$ of the local attractors

$$A_{class} = \overline{\overline{M_{2d}^u(p_0)} \backslash (M_{2d}^u(p_0) \cup \{p_0, p_1, p_2\})} \qquad (5)$$

When $A_{class}$ is nonempty, then a closed subset $A$ of $A_{class}$ is naturally called a non-trivial local attractor, if there is an integral curve $l$ dense on it, i.e., $\bar{l} = A$.

The numerical results will show that it is quite possible that there are at least two different local attractors in the same $A_{class}$.

The expression (5) means that the attractor class $A_{class}$ is constructed by the following steps: (i) make the closure of $M_{2d}^u(p_0)$, (ii) take $M_{2d}^u(p_0)$ and the three equilibrium points out of the closure, (iii) make the closure of the remained part of the closure previously obtained. Then the new closure is just the class $A_{class}$.

The closure operation in (5) guarantees the local attractiveness of these local attractors since, for each local attractor in this class, there is at least one integral on $M_{2d}^u(p_0)$ approaching to it.

It is believed that the system (2) is not integrable with quadrature, and that it is impossible to represent the integral curves, the manifolds mentioned above with the known functions.

Fortunately, based on the developments of the computer and numerical calculation techniques, it is now possible to study the system (2) numerically. First of all, it is possible to study the behaviors of these manifolds in a bounded region. To do so, just study the integral curves, which start at $t = 0$ from a very small neighborhood of a given equilibrium point. The half part of these curves, which correspond to the positive time $t$, may expand the unstable two-dimensional manifold, or may be treated as the approximations of the unstable one-dimensional manifold; and the other half part of the integral curves, which corresponds to the negative time $t$, may approximately expand the stable two-dimensional manifold, or may be treated as the approximations of the unstable one-dimensional manifold.

## § 2. Some Numerical Results

1. The pseudo gap of $M_{2d}^s(p_i)$.

If $0 \leq b < 3.1$, it seems that all of the stable manifolds and unstable manifolds of the three equilibrium points are unbounded, especially $M_{2d}^u(p_0)$ is unbounded.

Figures 1–4 show concretely these boundless manifolds for $b = 0.3$. So, it looks impossible to found out a bounded attractor in this parameter region.

For $0 \leq b < 0.202$, the branch $M^u_{1d}(p_1)_1$ will finally approach $M^u_{1d}(p_1)_2$; and symmetrically, the branch $M^u_{1d}(p_2)_1$ will finally approach $M^u_{1d}(p_2)_2$. From two different viewing angles, Figures 5 and 6 show how the branch $M^u_{1d}(p_2)_1$ will finally approach $M^u_{1d}(p_2)_2$ in the case $b = 0.2$.

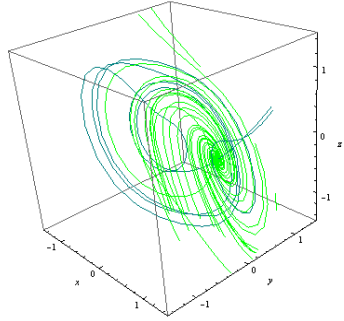
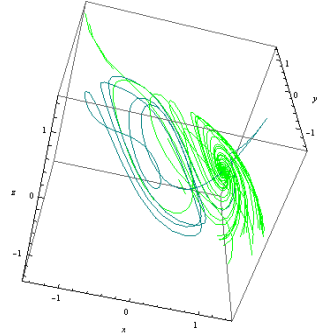

Figure 5. $M^u_{1d}(p_2)_1, M^u_{1d}(p_2)_2$ and $M^u_{2d}(p_2)$    Figure 6. $M^u_{1d}(p_2)_1, M^u_{1d}(p_2)_2$ and $M^u_{2d}(p_2)$

In the case that $M^u_{1d}(p_2)_1$ approaches $M^u_{1d}(p_2)_2$, the two-dimensional manifold $M^s_{2d}(p_2)$ seems to have a superficial "gap" to allow $M^u_{1d}(p_2)_1$ finally to pass $M^s_{2d}(p_2)$ through the "gap" to approach $M^u_{1d}(p_2)_2$. However, the "gap" does not mean that $M^s_{2d}(p_2)$ is really broken or fractured there; it only means that the manifold $M^s_{2d}(p_2)$ is strongly tortured and curved there to allow $M^u_{1d}(p_2)_1$ to pass superficially through the "gap" to approach , but not to intersect with, $M^u_{1d}(p_2)_2$. Symmetrically, the superficial "gap" of $M^s_{2d}(p_1)$ may appear when $M^u_{1d}(p_1)_1$ approaches $M^u_{1d}(p_1)_2$.

By similar reasoning, the superficial "gap" of $M^s_{2d}(p_1)$ may appear when $M^u_{1d}(p_2)_1$ approaches $M^u_{1d}(p_1)_2$, and the superficial "gap" of $M^s_{2d}(p_2)$ may appear when $M^u_{1d}(p_1)_1$ approaches $M^u_{1d}(p_2)_2$. In fact, this situation does exist when $b$ is near $0.3$.

When the parameter $b = 0.3$, $M^u_{1d}(p_2)_1$ finally will approach $M^u_{1d}(p_1)_2$, and $M^u_{1d}(p_1)_1$ finally will approach $M^u_{1d}(p_2)_2$. Figure 7 and 8 show this phenomenon.

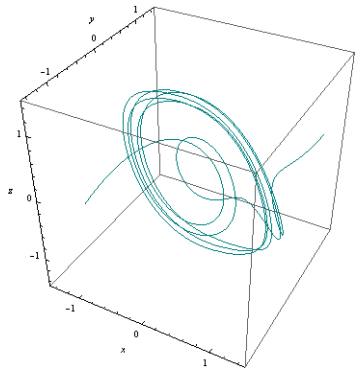
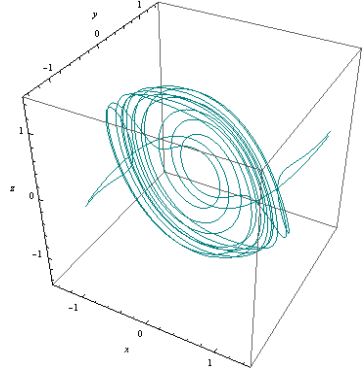

Figure 7. $M^u_{1d}(p_1)$    Figure 8. $M^u_{1d}(p_1)$ and $M^u_{1d}(p_2)$

The above-mentioned facts show that there must exist some particular parameter(s) $b^*$ ($0.202 < b^* < 0.3$), such that the branch $M_{1d}^u(p_1)_1$ (or $M_{1d}^u(p_2)_1$) becomes a homoclinic orbit from $p_1$ to $p_1$ (or from $p_2$ to $p_2$).

Because a large quantity of very fine numerical calculations are needed to represent the complex, tortured structure of $M_{2d}^s(p_1)$ and $M_{2d}^s(p_2)$ near the "gap", it is impractical for this paper to include that work.

Clearly, if the integral curves, which are located initially in the region between $M_{2d}^s(p_1)$ and $M_{2d}^s(p_2)$, approach infinity, they must pass one of the manifolds $M_{2d}^s(p_1)$ and $M_{2d}^s(p_2)$ through the corresponding "gap".

## 2. The Faint Attractor.

When $b = 0.311, 0.312$, the closure of $M_{2d}^u(p_0)$ seems to be a bounded set constrained by both $M_{2d}^s(p_1)$ and $M_{2d}^s(p_2)$. So by the expression (5), the system (2) may have a class $A_{class}$ of local attractors.

The numerical result shows that $A_{class}$ may be the only local attractor, though it has not been proven strictly.

The following figures are all based on the numerical results in the case $b = 0.312$.

Figures 9 and 10 show $M_{2d}^u(p_0)$ expressed through part of two integrals on it from two different viewing angles. And Figures 11 and 12 show the relations between $M_{2d}^u(p_0)$, $M_{2d}^s(p_1)$ and $M_{2d}^s(p_2)$ through two different viewing angles.

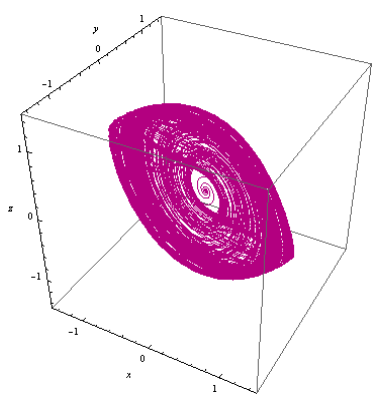
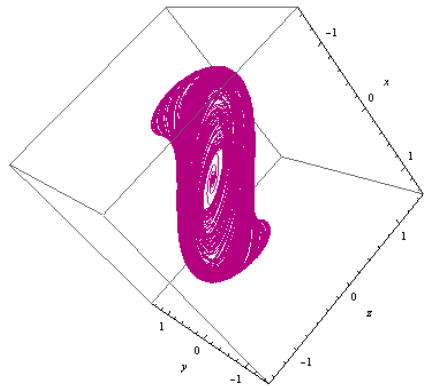

Figure 9. $M_{2d}^u(p_0)$    Figure 10. $M_{2d}^u(p_0)$

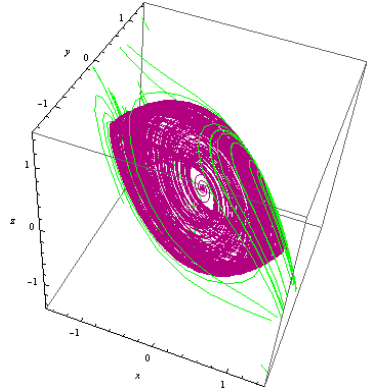
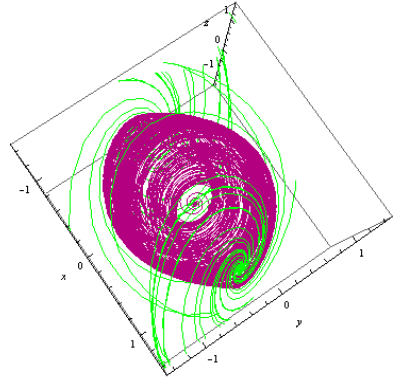

Figure 11. $M_{2d}^u(p_0)$, $M_{2d}^s(p_1)$ and $M_{2d}^s(p_2)$ 　　Figure 12. $M_{2d}^u(p_0)$, $M_{2d}^s(p_1)$ and $M_{2d}^s(p_2)$

Figures 13 and 14 show the set class $A_{class}$ from two different viewing angles.

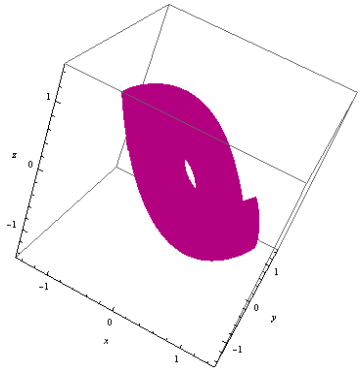
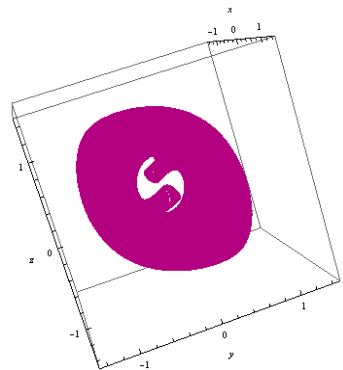

Figure 13. Set class $A_{class}$ 　　　　　　　　Figure 14. Set class $A_{class}$

The most interesting fact is that, in these cases, both $M_{1d}^u(p_1)_1$ and $M_{1d}^u(p_2)_1$ first approach $A_{class}$ and browse or visit almost over the whole $A_{class}$ in a very short distance, then finally will approach $M_{1d}^u(p_1)_2$ and $M_{1d}^u(p_2)_2$ respectively.

Figure 15 shows how $M_{1d}^u(p_1)_1$ begins to visit $A_{class}$, and Figure 16 shows how $M_{1d}^u(p_1)_1$ finally approaches $M_{1d}^u(p_1)_2$ after its browsing almost the whole part of $A_{class}$.

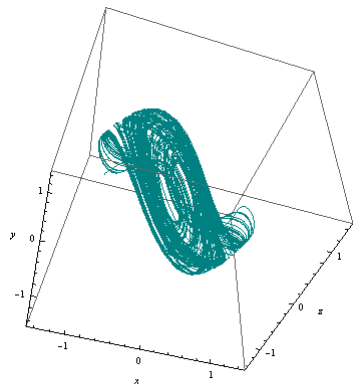
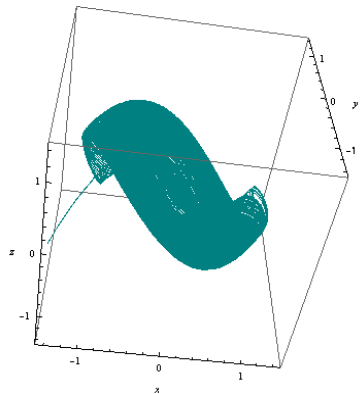

Figure 15. $M_{1d}^u(p_1)_1$ begins to visit $A_{class}$ 　　Figure 16. $M_{1d}^u(p_1)_1$ finally approaches $M_{1d}^u(p_1)_2$ after its browsing almost over the whole part of $A_{class}$

Figures 17 and 18 show how $M_{1d}^u(p_1)_1$ closes to $A_{class}$, respectively, from two different viewing angles, such that $M_{1d}^u(p_1)_1$ is almost merged with $A_{class}$.

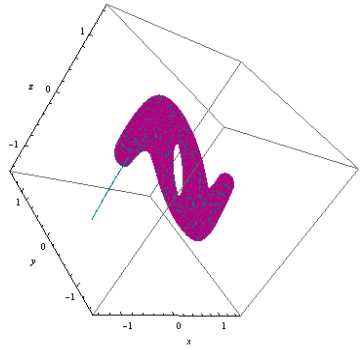 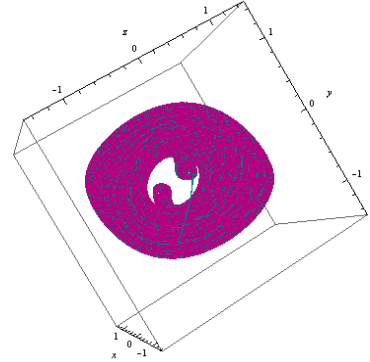

Figure 17.  $M_{1d}^u(p_1)_1$ and $A_{class}$        Figure 18.  $M_{1d}^u(p_1)_1$ and $A_{class}$

By the previous discussion, there must exist a "gap" on each of $M_{2d}^s(p_1)$ and $M_{2d}^s(p_2)$, so that there must exist a large quantity of integral curves located initially between $M_{2d}^s(p_1)$ and $M_{2d}^s(p_2)$ will finally approach infinity by passing through one of the two "gaps".

Because both $M_{1d}^u(p_1)_1$ and $M_{1d}^u(p_2)_1$ browse or visit almost over the whole $A_{class}$ in such a very short distance that they are almost merged with $A_{class}$, the set $A_{class}$ has only a very faint attraction. In this case, the local attractor in $A_{class}$ is called a **faint attractor**.

Most of integral curves that start initially from the region between $M_{2d}^s(p_1)$ and $M_{2d}^s(p_2)$ will first browse more or less part of $A_{class}$, and then go through one gap and approach infinity.

Figure 19 shows an integral curve starting from the initial point $(1.45, 0, -1.5)$ and finally approaching infinity through the gap of $M_{2d}^s(p_2)$. Figure 20 shows an integral curve starting from the initial point $(0.99, 0, 0)$ and finally approaching infinity through the gap of $M_{2d}^s(p_1)$.

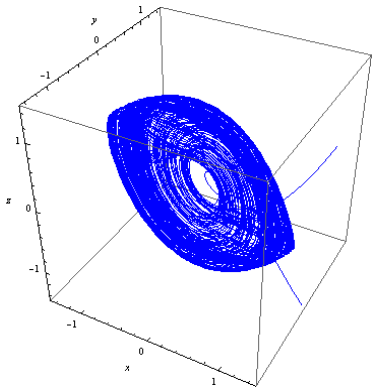 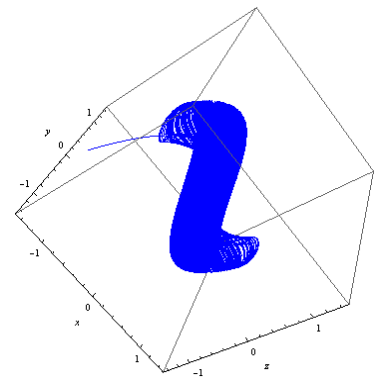

Figure 19. Integral curve starting from (1.45, 0, -1.5)    Figure 20. Integral curve starting from (0.99, 0, 0)

If an integral curve starts from an initial point outside the region between $M_{2d}^s(p_1)$ and $M_{2d}^s(p_2)$, it will directly approach infinity without browsing the faint $A_{class}$.

Figure 21 shows an integral curve starting from the point $(1.5, -1.5, 0)$, which is outside the region between $M_{2d}^s(p_1)$ and $M_{2d}^s(p_2)$, never browses the faint attractor.

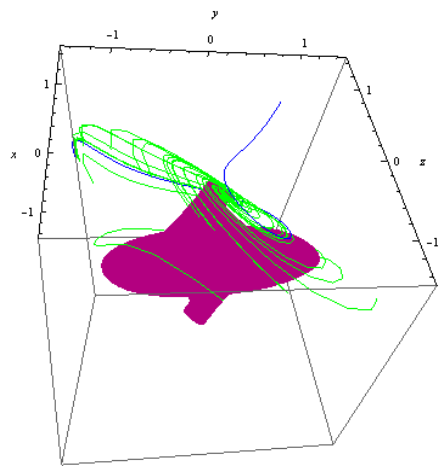

Figure 21. An integral curve never visiting the local attractor

3. Non-trivial Local Attractor.

If $b = 0.313, 0.314, 0.315, 0.316, 0.317$, then the closure of $M_{2d}^u(p_0)$ is clearly bounded, and the set $A_{class}$ is a local attractor. Besides, both $M_{1d}^u(p_1)_1$ and $M_{1d}^u(p_2)_1$ are bounded with $A_{class}$ as their limit set, i.e.,

$$A_{class} \subset \overline{M_{1d}^u(p_i)_1} \qquad i = 1, 2. \qquad (6)$$

Therefore, the superficial gaps of $M_{2d}^s(p_1)$ and $M_{2d}^s(p_2)$, which appeared in the previous cases, disappear now. All of the integral curves, which start at initial points located in the region between $M_{2d}^s(p_1)$ and $M_{2d}^s(p_2)$, will be attracted to the attractor $A_{class}$. The integral curves outside the region between $M_{2d}^s(p_1)$ and $M_{2d}^s(p_2)$ are unbounded.

The following figures are based on the numerical calculation in the case $b = 3.17$.

Figure 22 shows $A_{class}$, $M_{2d}^s(p_1)$ and $M_{2d}^s(p_2)$. Figure 23 shows $M_{1d}^u(p_2)_1$, $M_{2d}^s(p_1)$ and $M_{2d}^s(p_2)$. Figure 24 shows how $M_{1d}^u(p_2)_1$ approaches $A_{class}$. Figure 25 shows how an integral curve that starts from $(1.2, 0, -1.5)$, approaches $A_{class}$

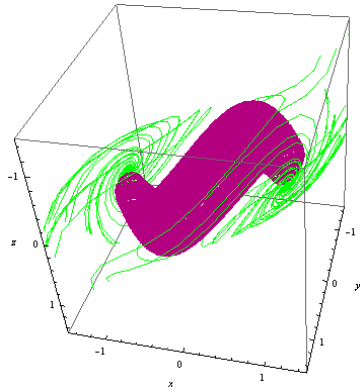
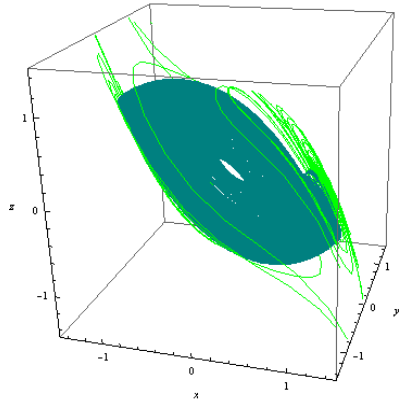

Figure 22. $A_{class}$, $M_{2d}^s(p_1)$ and $M_{2d}^s(p_2)$      Figure 23. $M_{1d}^u(p_2)_1$, $M_{2d}^s(p_1)$ and $M_{2d}^s(p_2)$

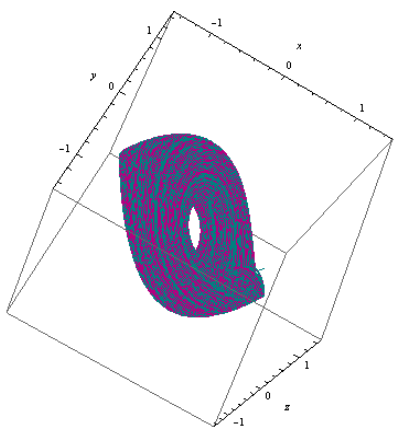
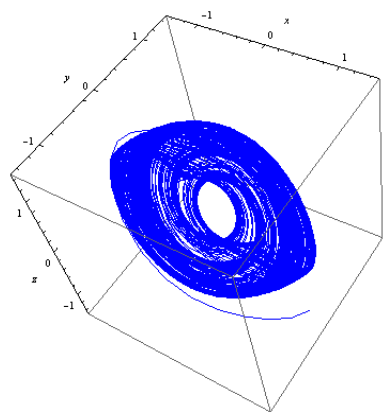

Figure 24. $A_{class}$ and $M_{1d}^u(p_2)_1$      Figure 24. An integral curve starting from $(1.2, 0, -1.5)$

The numerical results show that $A_{class}$ has a complex (maybe fractal) structure. Figure 26 shows the cross section of $A_{class}$ at the plane $x = 0$

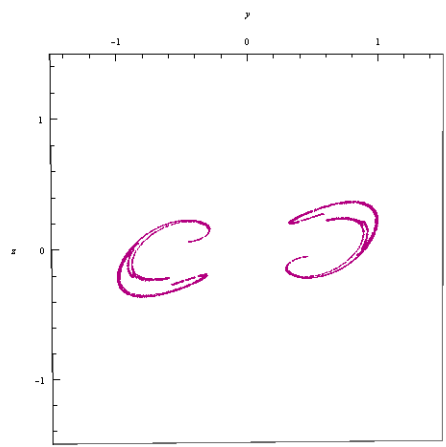

Figure 26. The cross section of $A_{class}$ at the plane $x = 0$

## 4. Twin Spatial Limit Closed Orbits and Bifurcation of Multiple Rotation Closed orbits.

If $b = 0.3175, 0.318, 0.319, 0.32$, the set $A_{class}$ seems to be formed obviously with a pair of spatial limit closed orbits, and they are symmetric about the origin.

Figures 27, 28, 29 and 30 show, respectively from two different viewing angles, the two closed spatial orbits and their relation, where the parameter $b = 0.318$.

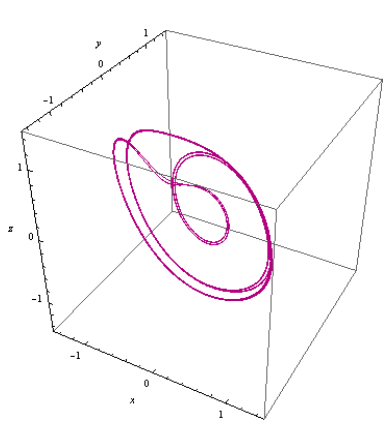
Figure 27. A spatial limit closed orbit

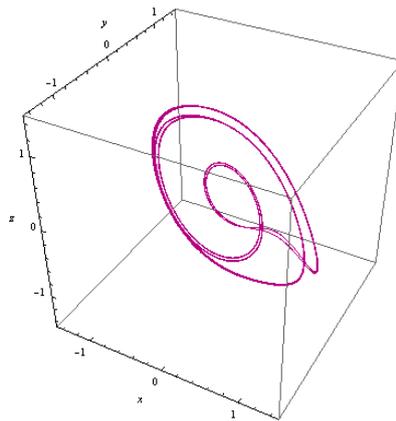
Figure 28. Another spatial limit closed orbit

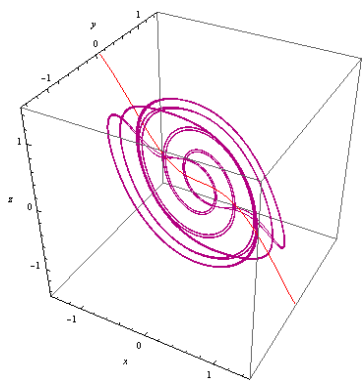
Figure 29. The twin spatial limit closed orbits and $M_{1d}^s(p_0)$

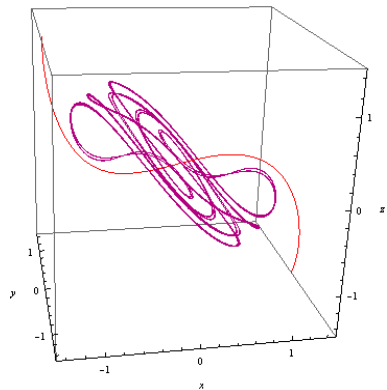
Figure 30. The twin spatial limit closed orbits and $M_{1d}^s(p_0)$

Especially, in Figures 29 and 30, the one-dimensional stable manifold $M_{1d}^s(p_0)$ is shown for the study of the topological character of the twin limit closed orbits.

In the case $b = 0.32$, Figures 31, 32, 33 and 34 show the twin spatial limit closed orbits in the same way as Figures 27, 28, 29 and 30.

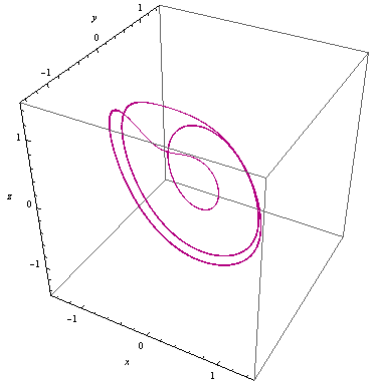 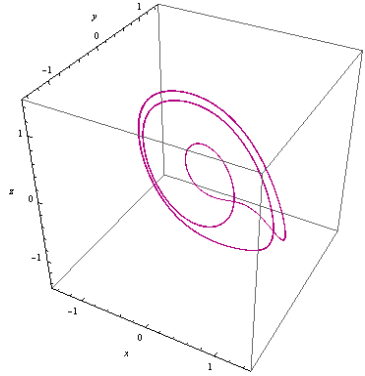

Figure 31. $b = 0.32$          Figure 32. $b = 0.32$

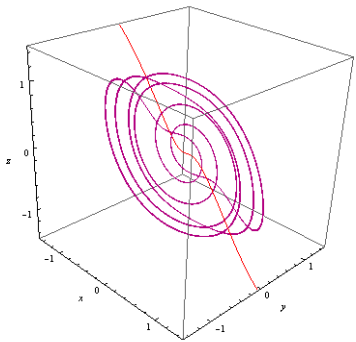 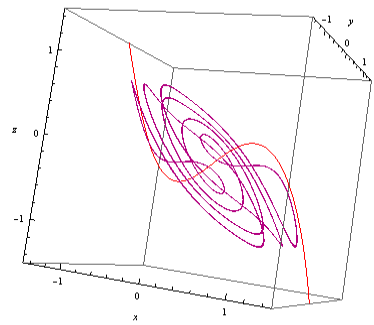

Figure 33. Twin closed orbits, $b = 0.32$     Figure 34. Twin closed orbits, $b = 0.32$

    Comparing the two groups of figures, it is easy to see that a moving point along a limit closed orbit must rotate around the one-dimensional stable manifold $M_{1d}^{s}(p_0)$ for a fixed number of times to meet its initial position. The number of the rotations is called the rotation number of the limit closed orbits.

    The rotation number is 6 in the case $b = 3.18$, and the rotation number is 3 in the case $b = 0.32$. In the case $b = 0.3175$, though it is difficult to see the exact rotation number clearly from Figure 35, however it is easy to see that the rotation numbers of the spatial closed orbit must be a multiple of 3, say $3 \times 2^{n^{*}}$.

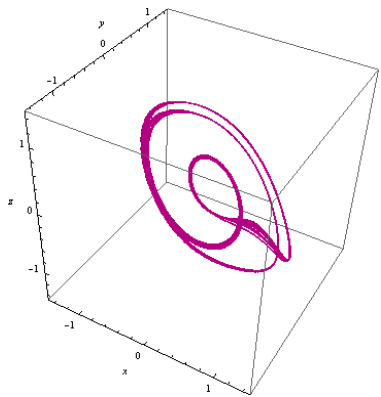

Figure 35. A spatial limit closed orbit with rotation number $3 \times 2^{n^{*}}$, $b = 0.3175$.

These facts show also that the rotation number of a spatial limit closed orbit decreases to three as the parameter $b$ increases from $0.3175$ to $0.32$.

Since the one-dimensional stable manifold $M_{1d}^s(p_0)$ cannot intersect with the spatial limit closed orbit, this paper proposes that there must exist a series of critical parameters $b_i^*$ $(i = 0, 1, 2, \dots, n^* - 1)$

$$0.315 < b_{n^*-1}^* < b_{n^*-2}^* < \cdots < b_0^* < 0.32,$$

such that, if $0.317 \le b < b_{n^*-1}^*$, then the rotation number of the spatial limit closed orbit is equal to $3 \times 2^{n^*}$; if $b_k^* < b < b_{k-1}^*$, then the rotation number is equal to $3 \times 2^k$ $(k = 1, \dots, n^* - 1)$; and if $b_0^* < b \le 0.32$, then the rotation number is equal to $3$.

The present numerical test shows that, $b_0^* \approx 0.3184$.

At these critical parameters, the system (2) also has a kind of structural instability on the rotation number. They may be called the **bifurcation points of the rotation number of the spatial limit closed orbits.**

The attracting region for these twin limit curves is naturally the region between $M_{2d}^s(p_1)$ and $M_{2d}^s(p_2)$. However, since the twin limit orbits are symmetric about the origin, and the attracting region must be separated into two symmetric parts such that each part corresponds to its own limit curve, so there must exist a nonempty boundary $B$ between the two symmetric regions. However, this paper has not probed what this boundary looks like.

5. New Parameter Region of Non-trivial Local Attractor

If $b = 0.33$, the situation is just like the case $b = 0.313, 0.314, 0.315, 0.316, 0.317$, the set $A_{class}$ is again a non-trivial local attractor with complex structure. Figures 36 and 37 show the set $A_{class}$ from two different viewing angles.

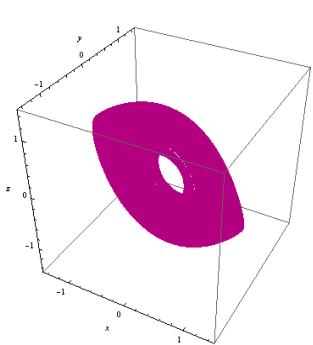 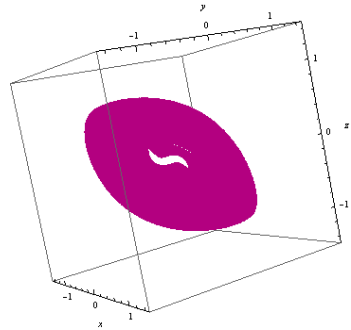

Figure 36. $A_{class}$, $b = 0.33$     Figure 37. $A_{class}$, $b = 0.33$

## 6. Single Spatial Limit Closed Orbit

If $b = \frac{1}{3}$, then the set $A_{class}$ seems to be formed with a single spatial limit closed orbit with a rotation number of 13.  See Figure 38, 39.

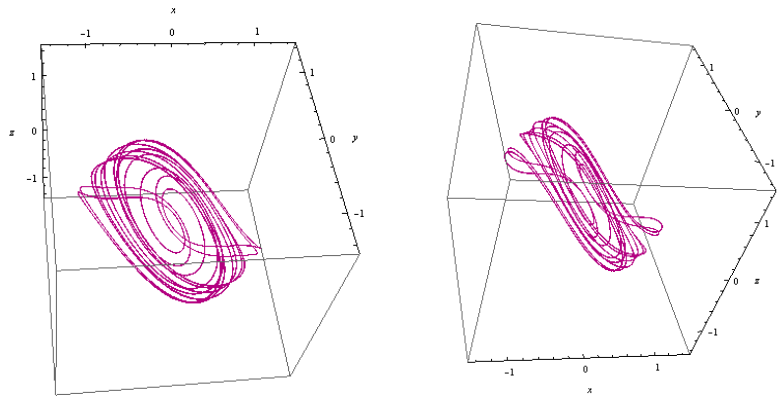

Figure 38, 39. Limit Closed Orbit of Rotation number 13.

The singleness is proposed here by that this limit curve seems symmetrical about the origin, and that different numerical calculations give the same limit orbit.

It is to be emphasized here that the rotation number of the limit orbit is slightly sensitive to the parameter $b$ around $b = \frac{1}{3}$.  For instance, when $b = 0.3332$, the rotation number of the limit curve is clearly larger than 13.  See Figure 40.

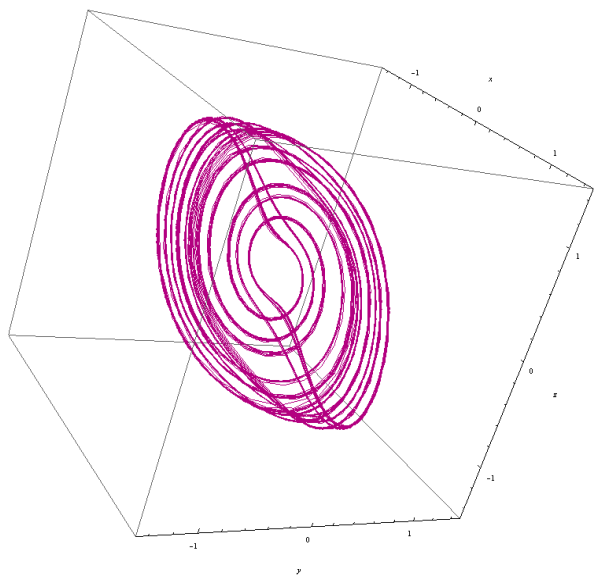

Figure 40. Limit closed orbit with larger rotation number, $b = 0.3332$

## 7. Newer Parameter Region of Non-trivial Local Attractor

If $b = 0.34$, the set $A_{class}$ is still like the one in the case $b = 0.33$.  See Figure

41.

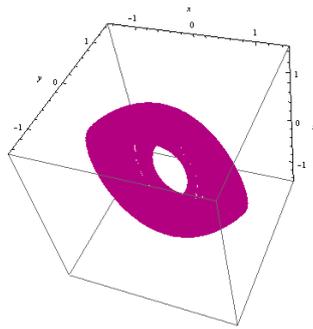

Figure 41. $A_{class}$, $b = 0.34$

8. Twin Spatial Limit Cycle

If $b = 0.4$, then $A_{class}$ is formed with a pair of symmetrical spatial limit cycles, and the rotation number for each cycle is 1. See Figure 42.

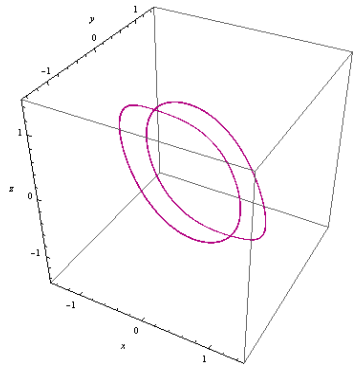

Figure 42. The pair of spatial limit cycles, $b = 0.4$

9. Single Spatial Limit Cycle

If $0.49 \leq b < 1$, then $A_{class}$ is just a single spatial limit cycle with rotation number of 1.

Figures 43 and 44 shows the spatial limit cycle from two different viewing angles, where $b = 0.5$.

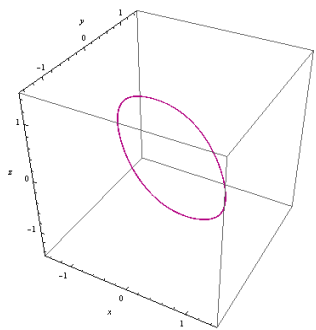 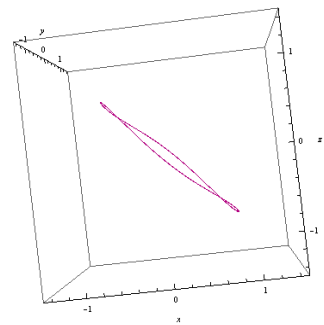

Figures 43. Spatial limit cycle, $b = 0.5$  Figures 44. Spatial limit cycle, $b = 0.5$

When the parameter $b$ is increasing from $0.5$ to $1$, the size of the limit cycle gets smaller and smaller. Figure 45 shows the limit cycle in the case $b = 0.99$. Figure 46 shows how an integral curve converges to the limit cycle.

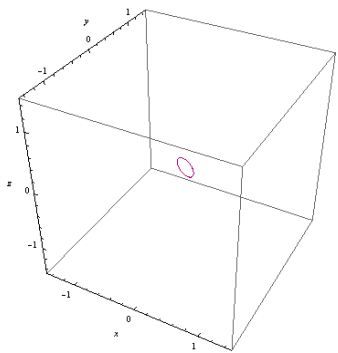
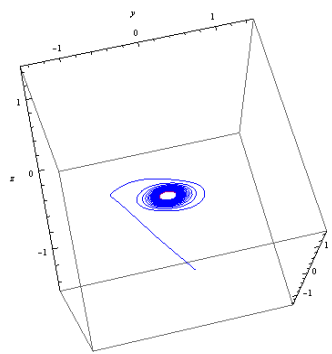

Figure 45. Limit cycle, $b = 0.99$    Figure 46. An integral curve converges to the limit cycle

There must exist such a bifurcation point of the parameter $b$, that around this parameter, the twin spatial limit cycles merge to a single limit cycle. The numerical tests show that this parameter is near and perhaps a little less than $b = 0.49$.

The parameter $b = 1$ is just the Hopf bifurcation point such that a small limit cycle appears when the parameter $b$ is little less than $1$.

## § 3.  More Discussion

Obviously, the numerical results in §2 cannot be treated as strictly mathematical conclusions.  However, they should be treated seriously, since they do provide us a series of new phenomena in the study of dynamical system theory based on a concrete ODE system.

Usually, it is relatively easy to consider a theory based on some general assumptions about abstract dynamical systems, but it is very difficult to give a concrete ODE to realize the ideal theory.  In fact, a concrete ODE system with interesting properties is an important source in the research. The history of the study on Lorenz equation is the best example.

The results of this paper show that the ODE system (2) is also an important source in the study of the dynamical system theory.

Recently, Prof. K. Kassner suggested M. Bradley change his ODE (1) into

$$y''' + y' = y(y^2 - a^2) \qquad (1')$$

I found that this leads the system (2) to be changed into

$$\begin{cases} \frac{dx}{dt} = P(x,y,z) = y \\ \frac{dy}{dt} = Q(x,y,z) = z \\ \frac{dz}{dx} = R(x,y,z) = a^2 x - x^3 - y - b z \end{cases} \qquad (2')$$

and that the new system (2') is more like the Lorenz equation (3). For instance, if $a = 1, b = 0.86$, the new system (2') has an attractor similar to the Lorenz attractor. See Figure 47.

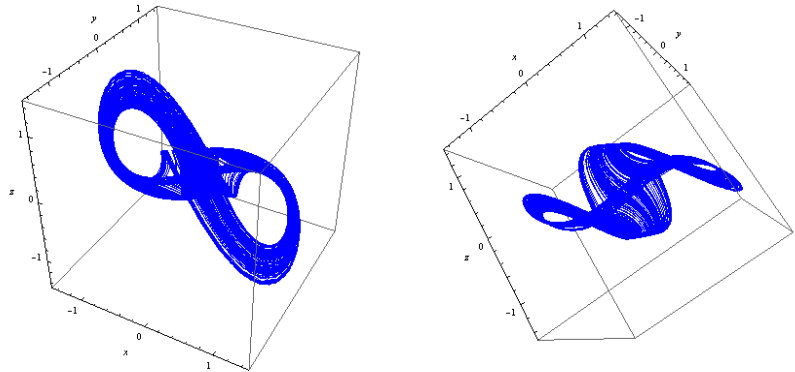

Figure 47. An attractor similar to Lorenz attractor.

In the discussion on the problem of Prof. R. Mark Bladley, the author has noticed that Dr. Julien Clinton Sprott had studied the system (2) numerically in the particular case, $a = 1, b = 0.7$, in his book[2]. He mentioned that the system (2) has a strange attractor with Lyapunov exponents $0.1380, 0, -0.8380$ for initial conditions of $(0,0,0.1)$. However, this result is in some doubt. Because the attractor of the system (2) should be a spatial limit cycle, not a strange attractor. Maybe the system studied by Dr. Julien Clinton Sprott is not (2), but just (2'), for the system (2') has just an attractor similar to the one shown in Figure 47.

Clearly, there must be more interesting phenomena and relevant theories about the system (2) and (2') to be uncovered or to be developed. I propose that these systems be explored more deeply in the future.

### Acknowledgements

Henry Ruddle, my friend, has read carefully the whole draft of this paper. He has revised some written errors and gave me a series of suggestions for improvements on the English, and on the expression of my idea. That is really a great help to me. The present paper has absorbed all of his suggestions. To him, I show my hearty thanks!

References

[1] John Guckenheimer and Philip Holmes, Nonlinear Oscillations, Dynamical Systems, and Bifurcations of Vector Fields, Springer-Verlag, New York, 1983

[2] Julien Clinton Sprott, Chaos and Time-series Analysis, Oxford University Press, 2003

**Postscript**

After this paper has been completed, the author noticed that the system (2) is a particular Silnikov equation

$$\begin{cases} \frac{dx}{dt} = P(x,y,z) = y \\ \frac{dy}{dt} = Q(x,y,z) = z \\ \frac{dz}{dx} = R(x,y,z) = \alpha(\delta x^3 - \gamma x^2 - x) - y - \beta z \end{cases} \quad (7)$$

i.e., if let $\alpha = a^2, \delta = \frac{1}{a^2}, \gamma = 0, \beta = b$, then the system (7) becomes the system (2).

From the literatures that I can find, the system (7) has been mentioned only on the webpage: SIMLIB Example: Silnikov equation,

http://www.fit.vutbr.cz/~peringer/SIMLIB/examples/silnikov.html

In this webpage, the system (7) has been studied numerically when $\alpha$ is negative, concretely, when $\alpha = -0.65, \beta = 0.4, \gamma = 0, \delta = 1$, and a strange attractor has been found. See Figure 48.

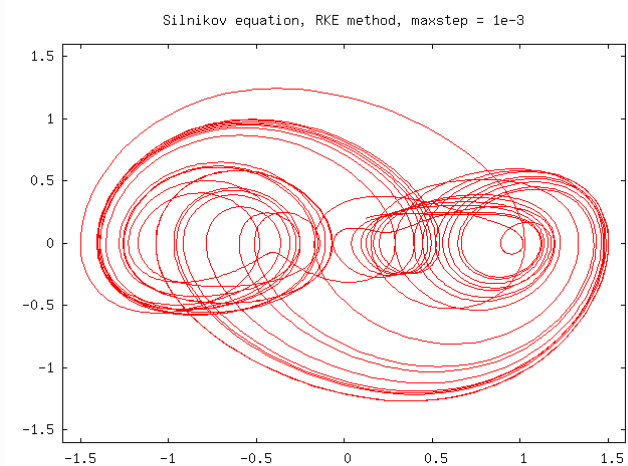

Figure 48. Attractor of Sinikov equation

It is not difficult to see that this particular Silnikov equation is almost equivalent to the system (2'). Maybe, the result of Dr. Julien Clinton Sprott mentioned above belongs to the case of negative α also.

Anyway, the present research paper shows the importance of the study of the Silnikov equation (7), especially for the positive α case.